\documentstyle[12pt]{article}
\catcode`\@=11

\def\bn{\hbox{\it I\hskip -2pt N}}
\def\bz{\hbox{\it Z\hskip -4pt Z}}

\def\demo{\noindent{\bf Proof .}}
\def\section{\@startsection {section}{1}{\z@}{-3.5ex plus -1ex
minus-.2ex}{2.3ex plus .2ex}{\normalsize\bf}}
\newtheorem{Theorem}{Theorem}
\newtheorem{Lemma}{Lemma}
\newtheorem{Proposition}{Proposition}

\newcommand{\het}{H_{\rm et}}
\newcommand{\hc}{H_{\rm c}}

\begin{document}
\begin{center}
{\LARGE\bf \textsc{On toric varieties of high arithmetical rank}}
\end{center}
\vskip.5truecm
\begin{center}
{Margherita Barile\footnote{Partially supported by the Italian Ministry of Education, University and Research.}\\ Dipartimento di Matematica, Universit\`{a} di Bari,Via E. Orabona 4,\\70125 Bari, Italy}
\end{center}
\vskip1truecm
\noindent
{\bf Abstract} We describe a class of toric varieties in the $N$-dimensional affine space which are minimally defined by no less than $N-2$ binomial equations.\vskip0.5truecm
\noindent
\section*{Introduction} The {\it arithmetical rank} (ara) of an algebraic variety  is the minimum number of equations that are needed to define it set-theoretically. For every affine variety $V\subset K^N$ we have that codim\,$V\leq$\,ara\,$V\leq N$. This general upper bound was found by Eisenbud and Evans \cite{EE}. In particular cases a better upper bound can be obtained by direct computations based on Hilbert's Nullstellensatz: this was done for certain toric varieties in \cite{B1}, \cite{B2}, \cite{BL}, \cite{BMT1}, \cite{BMT2}. In all the examples treated there the arithmetical rank was close to the trivial lower bound, i.e., ara\,$V\leq$\,codim\,$V+1$.  In this paper we present a class of toric varieties whose arithmetical rank is close to the general upper bound, namely  ara\,$V\geq N-2$. For proving this result, of course, we need a more efficient lower bound: this is provided by \'etale cohomology. The same kind of tools was used in \cite{BL}. There they were applied for showing that the arithmetical rank of certain toric varieties of codimension 2 depends on the characteristic of the ground field, and that ara\,$V=$\,codim\,$V$ in exactly one positive characteristic. In the present paper, however, we study toric varieties of any codimension, and obtain the same lower bound in all characteristics. This lower bound turns out to be sharp if $N=3$ or $N=5$.

\section{The main theorem}
Let $K$ be an algebraically closed field, and let $n\geq2$ be an integer. Let ${\bf e}_1,\dots, {\bf e}_n$ be the standard basis of $\bz^n$. Set $N=2n-1$ and consider the following subset of $\bn ^N$:
$$T=\{{\bf e}_1,\dots, {\bf e}_{n-1}, d{\bf e}_n, a_1{\bf e}_1+{\bf e}_n,\dots, a_{n-1}{\bf e}_{n-1}+{\bf e}_n\},$$
\noindent 
where $d,a_1,\dots, a_{n-1}\in\bn^{\ast}$. We also suppose that $d$ is divisible by two distinct primes $p$ and $q$. With $T$ we can associate the variety $V$ admitting the following parametrization
$$V:\left\{
\begin{array}{rcl}
x_1&=&u_1\\
&\vdots&\\
x_{n-1}&=&u_{n-1}\\
x_n&=&u_n^d\\
y_1&=&u_1^{a_1}u_n\\
&\vdots&\\
y_{n-1}&=&u_{n-1}^{a_{n-1}}u_n
\end{array}\right.,$$
which is a toric variety of codimension $n-1$ in the affine space $K^N$. Our aim is to show the following 
\begin{Theorem}\label{main}
{\rm ara}\,$V\geq N-2$.
\end{Theorem}
This will be done in the next section.
\section{The lower bound}
We show that ara\,$V\geq N-2$ by means of the following criterion, which is based on \'etale cohomology and is cited from \cite{BS}, Lemma 3$^\prime$. 
\begin{Lemma}\label{Newstead}Let
$W\subset\tilde W$ be affine varieties. Let $d=\dim\tilde
W\setminus W$. If there are $s$ equations $F_1,\dots, F_s$ such
that $W=\tilde W\cap V(F_1,\dots,F_s)$, then 
$$\het^{d+i}(\tilde W\setminus W,{\bz}/r{\bz})=0\quad\mbox{ for all
}i\geq s$$ and for all $r\in{\bz}$ which are prime to {\rm char}\,$K$.
\end{Lemma}
Since $p$ and $q$ are distinct primes, we may assume that char\,$K\ne p$. Hence our claim will follow once we have shown 
\begin{Proposition}\label{pq} If {\rm char}\,$K\ne p$, then 
$$\het^{2N-3}(K^N\setminus V, \bz/p\bz)\neq0.$$
\end{Proposition}
\demo In the sequel $\het$  and $\hc$ will denote \'etale cohomology and \'etale cohomology with compact support with respect to the coefficient group $\bz/p\bz$: we shall omit the latter for the sake of simplicity. By Poincar\'e Duality (see \cite{M}, Theorem 14.7, p.~83) we have 
\begin{equation}\label{7}{\rm Hom}_{\bz}(\het^{2N-3}(K^N\setminus V),\bz/p\bz)\simeq\hc^3(K^N\setminus V).\end{equation}
\noindent
Moreover, we have a long exact sequence of cohomology with compact support
$$\cdots\rightarrow\hc^2(K^N)\rightarrow\hc^2(V)\rightarrow\hc^3(K^N\setminus V)\rightarrow\hc^3(K^N)\rightarrow\cdots,$$
\noindent
where $\hc^2(K^N)=\hc^3(K^N)=0$, since $N\geq3$ (see \cite{M}, Example 16.3, pp.~98--99). Hence
\begin{equation}\label{8} \hc^3(K^N\setminus V)\simeq\hc^2(V).\end{equation}
\noindent
By (\ref{7}) and (\ref{8}) it thus suffices to show that 
\begin{equation}\label{9} \hc^2(V)\neq0.\end{equation}
\noindent
On $K^n$ fix the coordinates $u_1,\dots, u_n$ and let $X$ be the subvariety of $K^n$ defined by $u_1=u_2=\cdots=u_{n-1}=0$. Then $X$ is a 1-dimensional affine space over $K$, on which we fix the coordinate $u_n$. 
Consider the surjective map
$$\phi:K^n\rightarrow V$$
$$(u_1,\dots, u_n)\mapsto(u_1,\dots, u_{n-1}, u_n^d, u_1^{a_1}u_n,\dots, u_{n-1}^{a_{n-1}}u_n)$$
\noindent
and the restriction map
$$\tilde\phi:K^n\setminus X\rightarrow V\setminus\phi(X),$$
\noindent
which is a bijective morphism of affine schemes.
For all $i=1,\dots,n-1$ let $V_i=\{(\bar x_1,\dots,\bar x_{n-1},\bar x_n,\bar y_1,\dots, \bar y_{n-1})\in V\mid\bar x_i\ne0\}$. These sets form an open cover of $V\setminus\phi(X)$, and $U_i=\phi^{-1}(V_i)=\{(u_1,\dots, u_n)\in K^n\mid u_i\ne0\}$.  Moreover, for all $i=1\,\dots, n-1$, the morphism
$$\psi_i:V_i\rightarrow U_i$$
$$(\bar x_1,\dots,\bar x_{n-1},\bar x_n,\bar y_1,\dots, \bar y_{n-1})\mapsto(\bar x_1,\dots,\bar x_{n-1}, \frac{\bar y_i}{\bar x_i^{a_i}})$$
\noindent is the inverse map of the restriction of $\tilde\phi$ to $U_i$. Hence $\tilde\phi$ is an isomorphism of affine schemes, so that it induces
an isomorphism of groups
$$\tilde\phi^{\ast}_i:\hc^i(V\setminus\phi(X))\simeq\hc^i(K^n\setminus X)$$
\noindent for all indices $i$. The restriction of $\phi$ to $X$   
$$\bar\phi:X\rightarrow\phi(X)$$
\noindent maps $u_n$ to $(0,\dots,0,u_n^d,\dots,0)$. Hence $\phi(X)$ is a 1-dimensional affine space and $\bar\phi$ induces multiplication by $d$ on the second cohomology group with compact support (see \cite{M}, Remark 24.2 (f), p.~135). Now, as it is well-known (see \cite{M}, Example 16.3, pp.~98--99),
$$\hc^i(\phi(X))\simeq\hc^i(X)\simeq\left\{\begin{array}{rl}\bz/p\bz&\mbox{ for }i=2\\
0&\mbox{else}\end{array}\right..
$$
\noindent
Since $p$ divides $d$, it follows that the induced maps 
$$\bar\phi^{\ast}_i:\hc^i(\phi(X))\longrightarrow\hc^i(X)$$
\noindent are all equal to the zero map. 
 Thus $\phi$ gives rise to the following morphism of acyclic complexes:
$$\begin{array}{cccccc}
&&\bz/p\bz&&&\\
&&\wr\vert&f&&\\
\hc^2(V)&\to&\hc^2(\phi(X))&\to&\hc^3(V\setminus\phi(X))\\
&&&&&\\
\downarrow&&\bar\phi^{\ast}_2\downarrow0&&\wr\vert\downarrow\tilde\phi^{\ast}_3\\
&&&&&\\
\hc^2(K^n)&\to&\hc^2(X)&\to&\hc^3(K^n\setminus X)\\
&&\wr\vert&g&&\\
&&\bz/p\bz&&&
\end{array}.$$
\noindent
Note that $\tilde\phi^{\ast}_3f=g0=0$, so that $\tilde\phi^{\ast}_3f$ is not injective. Since $\tilde\phi^{\ast}_3$ is an isomorphism, it follows that $f$ is not injective. This implies (\ref{9}) and completes the proof. 
\section{On the defining equations}
Finally we show that the lower bound established in Theorem \ref{main} is sharp, since it is attained by the arithmetical rank when $n=2$ or $n=3$.
\begin{Proposition} If $n=2$, $V$ is defined set-theoretically by $F=y_1^d-x_1^{a_1d}x_2$. If $n=3$, then $V$ is set-theoretically defined by the following three binomials:
$$F_1=y_1^d-x_1^{a_1d}x_3,\quad F_2=y_2^d-x_2^{a_2d}x_3,\quad G=y_1^{d-1}y_2-x_1^{a_1(d-1)}x_2^{a_2}x_3.$$
\end{Proposition}
\demo If $n=2$, the claim is clear, since the defining ideal of $V$ is the principal ideal generated by $F$. So suppose that $n=3$.  It is straightforward to check that for all ${\bf w}\in V$, $F_1({\bf w})=F_2({\bf w})=G({\bf w})=0$. Conversely, we have to prove that for every ${\bf w}\in K^5$ such that 
\begin{equation}\label{1}F_1({\bf w})=F_2({\bf w})=G({\bf w})=0,\end{equation}
\noindent we have that ${\bf w}\in V$.
Let ${\bf w}=(\bar x_1, \bar x_2, \bar x_3, \bar y_1, \bar y_2)\in K^5$ be a point fulfilling (\ref{1}). Set $u_i=\bar x_i$ for $i=1,2$. We show that, for a suitable choice of parameter $u_3$, we can write ${\bf w}=(u_1, u_2, u_3^d, u_1^{a_1}u_3, u_2^{a_2}u_3)$. This is certainly true if $\bar x_3=0$: in this case (\ref{1}) implies that 
$\bar y_1=\bar y_2=0$, and  $u_3=0$ is the required parameter. Suppose that $\bar x_3\neq0$.   Let $u_3$ be a $d$-th root of $\bar x_3$. By (\ref{1}) we have that, for $i=1,2$, 
$$\bar y_i^d=u_i^{a_id}u_3^d,$$
\noindent
which implies that 
\begin{equation}\label{3}\bar y_i=u_i^{a_i}u_3\omega_i\end{equation}
\noindent for some $d$-th root $\omega_i$ of 1. On the other hand, from (\ref{1}) we also deduce that 
\begin{equation}\label{4}\bar y_1^{d-1}\bar y_2=u_1^{a_1(d-1)}u_2^{a_2}u_3^d.\end{equation}
\noindent
Note that, by (\ref{1}), since $\bar x_3\ne 0$, for $i=1,2$, we have $\bar x_i=0$ if and only if $\bar y_i=0$. If $\bar x_1=0$, set $u_3'=u_3\omega_2$. Then $\bar x_3=u_3'^d$, and, in view of (\ref{3}), $\bar y_i=u_i^{a_i}u_3'$ for $i=1,2$. Hence $u_3'$ is the required parameter.  Similarly one can reason if $\bar x_2=0$. So assume that $\bar x_i\ne 0$ for $i=1,2$. 
Replacing (\ref{3}) on the left-hand side of (\ref{4}) gives
\begin{equation}\label{5}u_1^{a_1(d-1)}u_3^{d-1}\omega_1^{d-1}u_2^{a_2}u_3\omega_2=u_1^{a_1(d-1)}u_2^{a_2}u_3^d.\end{equation}
\noindent
Since $u_1, u_2, u_3$ are non zero, from this we deduce that
\begin{equation}\label{6}\omega_1^{d-1}\omega_2=1,\end{equation}
\noindent which implies that 
$$\omega_1=\omega_2,$$
\noindent 
 i.e., $\omega_1$ and $\omega_2$ are both equal to the same $d$-th root $\omega$ of 1. Thus $u_3'=u_3\omega$ is the required parameter.   
\par\medskip\noindent
The computation of the arithmetical rank for $n\geq 4$ remains an open question.

\end{document}